\documentclass{article}
\usepackage[cp1250]{inputenc}
\usepackage[vcentermath]{youngtab}
\usepackage{amstext,amsbsy,amsmath,amscd,amssymb, amsthm, array, xy}
\usepackage{mathrsfs, color}

\input xy
\xyoption{all}

\begin{document}
\title{Symplectic spinor valued forms  and invariant operators acting between them}
\author{Svatopluk Kr\'ysl \footnote{{\it E-mail address}: krysl@karlin.mff.cuni.cz}\\ {\it \small   Charles University, Sokolovsk\'a 83, 186 75, Prague 8, Czech Republic.}}

\maketitle

\noindent\centerline{\large\bf Abstract} Exterior differential forms with values in the (Kostant's)
symplectic spinor bundle on a manifold with a given metaplectic structure are decomposed into
invariant subspaces. Projections to these invariant subspaces of a covariant derivative associated
to a torsion-free symplectic connection are described.

 \section{Introduction}

While the spinor twisted de Rham sequence for orthogonal spin
structures is well understood from the point of view of
representation theory (see, e.g., Delanghe, Sommen, Sou\v{c}ek
\cite{DSS}), its symplectic analogue seems to be untouched till
present days. In   Riemannian geometry, a decomposition of
spinor twisted de Rham sequence (i.e., exterior differential forms
with values in   basic spinor bundles) into invariant
parts is well known. Suppose a principal connection on the frame bundle of orthogonal repers (of the tangent bundle) is given. It induces in a canonical way
a covariant derivative on differential forms with values in the basic spinor bundles.
In this case, it is known, which parts of the covariant
derivatives acting between the spinor bundle valued forms are zero if we
restrict it to an invariant part of the sequence. Namely, the
covariant derivative maps each invariant part only in at most three
invariant parts sitting in the next gradation (some degeneracies on
the ends of the sequence could be systematically described).
  In   symplectic geometry, the first question which naturally
  arises is, what are the spinors for a symplectic Lie algebra. This
  question was successfully answered by Bertram Kostant in
  \cite{Kostant2}. He offered a candidate for symplectic spinors. We
  will call these spinors {\it basic symplectic spinors} and denote their underlying vector spaces
 $\mathbb{S}_+$ and $\mathbb{S}_-$. They are
  analogous to the ordinary orthogonal spinors in at least two
  following ways.
\begin{itemize}
\item[]  First, they could
  be found in a {\it symmetric algebra} of an isotropic subspace, while the orthogonal spinors could be
  found in  an {\it exterior algebra} of certain isotropic subspace.

\item[]   Second, the highest weights of the basic symplectic spinors are also
  half integral like the highest weights of orthogonal spinors (both with respect to the usual basis of
  the dual of an appropriate Cartan subalgebra).
\end{itemize}

  Unlike the orthogonal spinors, the symplectic ones are of
  infinite dimension, thus not so easy to handle.

       The main results of this article are the  decomposition of the spinor twisted
     de Rham sequence in the symplectic case and a theorem, which says
     that the image of each covariant derivative (associated to a symplectic  torsion-free connection and restricted to an invariant subspace)
     lies in  at most three invariant subspaces, i.e., a similar theorem to
     that one, which is valid in Riemannian geometry.
To derive the first mentioned result, we need to decompose the
ordinary exterior forms into irreducible summands over the
symplectic Lie algebra - a procedure, which is well known. We also
need to know, how to decompose a tensor product of such irreducible
summand and the basic symplectic spinor module $\mathbb{S}_+$. This
was done by Britten, Lemire, Hooper in \cite{Britten} and by Britten,
Lemire in \cite{BH} even in a more general setting. To derive the
second result (the description of the image of the covariant
derivative), we need one more ingredient. In particular, we should know,
how to decompose a tensor product of the defining representation
$\mathbb{V}$ of the symplectic Lie algebra and each infinite-dimensional
representation of the symplectic Lie algebra, which is an
irreducible summand in the space $\bigwedge^i \mathbb{V} \otimes
\mathbb{S}.$ These irreducible summands belong to a broader class of
infinite dimensional modules over a symplectic algebra, so called
{\it higher symplectic  spinor modules}, which are also known as
harmonic spinor representations in the
literature. 
To describe the decomposition of the tensor product of the defining
representation and the mentioned irreducible summand, we shall use a
theorem which was derived by the author in \cite{Krysl}.

    Investigation of the decomposition of the
    twisted de Rham complex for metaplectic structures
    has been motivated by a search for symplectic analogues of an (orthogonal) Dirac
    operator and its generalizations, which naturally appear in the twisted
    de Rham sequence in the orthogonal setting. Namely, it is known that the
    Dirac,  Rarita-Schwinger and twistor operators could be found in the twisted
     de Rham sequence for an orthogonal spin structure.
     The symplectic Dirac operator has been found by B. Kostant, see
     \cite{Kostant2}, and and has been studied intensively by many authors, see,
     e.g., Habermann \cite{H}, Klein \cite{Klein} and
     Kadl\v{c}\'{a}kov\'{a} \cite{Lenka}.
  We will recover all 
of these operators (Dirac,  Rarita-Schwinger and twistor) in a more systematic way
by an investigation of invariant differential operators appearing in
the symplectic spinor twisted de Rham sequence for  metaplectic  structures.
 No definition of the  symplectic
Rarita-Schwinger operator within mathematics is known to the author.
In physical literature, there are some references to  symplectic
Majorana fields or symplectic  Rarita-Schwinger fields, see Reuter
\cite{Reuter} and Green, Hull \cite{Green}, in the context of
super-gravity of strings.

 In the algebraic  part of this article (part 2), some
basic and known facts on higher symplectic spinor modules and decomposition of the mentioned tensor
products are written (Lemma 1, Theorem 1, Theorem 2). Besides these theorems, the
first main result (the decomposition of the symplectic spinor twisted de Rham sequence) is described
(Lemma 2) together with the theorem on the decomposition of the tensor product of the
defining representation and a higher symplectic spinor module
(Theorem 3). In this part, an information on intersection
of $\mathfrak{g}$-modules is written (Lemma 3).
 Third part of this article is the  geometrical one. It contains a general
 lemma on an image of a covariant derivative (Lemma 4) and the second
 main result (Theorem 4), namely the characterization of the image of a
 covariant derivative associated to a torsion-free symplectic
 connection.

\section{Higher symplectic spinor modules}

Let $(\mathbb{V},\omega)$ be a complex symplectic space of complex
dimension $2l,$ $l \in \mathbb{N}.$ Let $G=Sp(\mathbb{V},\omega)
\simeq  Sp(2l,\mathbb{C})$ be a complex   symplectic group of
$(\mathbb{V},\omega)$ and
$\mathfrak{g}=\mathfrak{sp}(\mathbb{V},\omega)  \simeq
\mathfrak{sp}(2l,\mathbb{C})$ its Lie algebra.\footnote{Different
choices of the symplectic form lead to isomorphic symplectic groups
and algebras.} Consider a Cartan subalgebra $\mathfrak{h}$ of the
symplectic Lie algebra is given together with a choice of positive roots
$\Phi^+$ of the system of all roots   $\Phi.$ The set of fundamental
weights $\{\varpi_i\}_{i=1}^l$ is then uniquely determined. For
later use, we shall need an orthogonal basis (with respect to the Killing
form on $\mathfrak{g}$), $\{\epsilon_i\}_{i=1}^l,$ for which
$\varpi_i=\sum_{j=1}^i\epsilon_j$ for $i=1,\ldots, l.$

For $\lambda \in \mathfrak{h}^*,$ let $L(\lambda)$ be the (up
to a $\mathfrak{g}$-isomorphism uniquely defined) irreducible module
with the highest  weight  $\lambda.$ If $\lambda$ happens to be
integral and dominant (wr. to the choice $(\mathfrak{h},\Phi^+)$),
i.e., $L(\lambda)$ is finite dimensional, we shall write
$F(\lambda)$ instead of $L(\lambda).$ Let $L$ be an arbitrary (finite or infinite dimensional) weight module over a complex simple Lie algebra.
 We call $L$ {\it module with bounded multiplicities}, if there is a $k\in \mathbb{N}_0,$ such that for each $\mu \in \mathfrak{h}^*,$ $\mbox{dim}L_{\mu}\leq k,$ where $L_{\mu}$ is the weight space of weight $\mu.$

Let us introduce the following set of weights
$$ A:=\{\lambda=\sum_{i=1}^{l}\lambda_i\varpi_i| \lambda_{l-1}+2\lambda_l + 3 >0,
\lambda_i \in \mathbb{N}_0, i=1,\ldots, l-1, \lambda_l \in \mathbb{Z}+\frac{1}{2}\}$$

{\bf Definition 1:} For a weight $\lambda \in
A,$  we call the  module $L(\lambda)$ {\it higher symplectic spinor module}. We shall denote the module
$L(-\frac{1}{2}\varpi_{l})$ by $\mathbb{S}_+$ or simply by
$\mathbb{S}$ and the module $L(\varpi_{l-1}-\frac{3}{2}\varpi_l)$ by
$\mathbb{S}_-.$ We shall call these two representations {\it basic
symplectic spinor modules}.

The next theorem says that the class of higher symplectic spinor
modules is quite natural and broad in a sense.

{\bf Theorem 1:} Let $\mathfrak{g}\simeq \mathfrak{sp}(2l,\mathbb{C})$ and $\lambda \in \mathfrak{h}^*$. Then the following are equivalent
\begin{itemize}
\item[1.)] $L(\lambda)$ is a module with bounded multiplicities
\item[2.)] $L(\lambda)$ is a direct summand in the completely reducible tensor product $\mathbb{S}\otimes F(\nu)$ for some integral dominant $\nu \in \mathfrak{h}^*$
\item[3.)] $\lambda \in A$
\end{itemize}
{\it Proof.} See Britten, Hooper, Lemire \cite{Britten} and Britten, Lemire \cite{BH}. $\Box$

In this paper, we shall first study the irreducible decomposition of
spaces $\bigwedge^{i} \mathbb{V}^*\otimes \mathbb{S}$ for
$i=0,\ldots, 2l.$  To do it, we need to decompose the wedge powers
$\bigwedge^i \mathbb{V}$ into irreducible modules. In the symplectic case (contrary to
the  orthogonal  one),  the wedge powers are not irreducible
generically. This decomposition
is well known and we state it as lemma 1.

{\bf Lemma 1:} Let $\mathbb{V}$ be the $2l$ dimensional defining representation of the symplectic Lie algebra $\mathfrak{sp}(\mathbb{V}, \omega),$ then
$$\bigwedge^{i}\mathbb{V} \simeq \bigoplus_{p=0}^{[i/2]}F(\varpi_{ i - 2p})$$ for $i=0,\ldots, l,$ 
where $[q]$ is the lower integral part of an element $q \in \mathbb{R}.$

{\it Proof.} See Goodman, Wallach \cite{GW}, pp. 237. $\Box$

In the next theorem, the decomposition of the tensor product of
an irreducible  finite dimensional $\mathfrak{sp}(\mathbb{V,\omega})$-module and
the basic symplectic spinor module  $\mathbb{S}$ is described.

{\bf Theorem 2:} Let $\mathfrak{g} \simeq \mathfrak{sp}(2l,\mathbb{C})$ and $\nu=\sum_{i=1}^l\nu_i\varpi_i \in \mathfrak{h}^*$ be an integral dominant weight for $(\mathfrak{h},\Phi^+).$ Let us define a
set $T_{\nu}:=\{\nu-\sum_{i=1}^ld_i\epsilon_i| d_i \in \mathbb{N}_0, \sum_{i=1}^l d_i \in 2\mathbb{Z}, 0\leq d_i \leq \nu_i, i=1,\ldots, l-1, 0\leq d_l \leq 2\nu_l+1\}.$ Then
$$F(\nu) \otimes \mathbb{S}=\bigoplus_{\kappa \in T_{\nu}}L(\kappa -\frac{1}{2}\varpi_l).$$

{\it Proof.} See Britten, Lemire \cite{BH}, theorem 1.2.$\Box$

We use the two last written claims to decompose the tensor products
$\bigwedge^i \mathbb{V}^* \otimes \mathbb{S},$ for $i=0,\ldots, l.$
We shall introduce the following convention.  For  a weight
$\sum_{i=1}^{l}\lambda_i \varpi_i,$ we shall write  $(\lambda_1
\lambda_2 \ldots \lambda_{l})$ briefly instead of
$L(\sum_{i=1}^l\lambda_i \varpi_i).$

{\bf Lemma 2:}
The following decompositions hold:

{\footnotesize For $2i+1\leq l-1,$
$$\bigwedge^{2i+1}\mathbb{\mathbb{V}^*}\otimes \mathbb{S} \simeq \bigwedge^{2l-2i-1}\mathbb{V}^*\otimes \mathbb{S} \simeq (10 \ldots
0-\frac{1}{2}) \oplus (0010 \ldots 0-\frac{1}{2})\oplus\ldots \oplus
(0 \ldots  0 1_{2i+1} 0 \ldots 0 -\frac{1}{2})$$
$$\oplus (0\ldots 0 1-\frac{3}{2})\oplus(010\ldots 01-\frac{3}{2})\oplus \ldots \oplus (0\ldots01_{2i}0\ldots01-\frac{3}{2})$$
For $2i\leq l-1$
$$\bigwedge^{2i }\mathbb{\mathbb{V}^*}\otimes \mathbb{S} \simeq \bigwedge^{2l-2i}\mathbb{V}^*\otimes \mathbb{S}\simeq (0 \ldots 0-\frac{1}{2}) \oplus (010 \ldots 0-\frac{1}{2})\oplus\ldots \oplus (0 \ldots  0 1_{2i} 0 \ldots 0 -\frac{1}{2})$$
$$\oplus (10\ldots  01-\frac{3}{2})\oplus(0010\ldots 01-\frac{3}{2})\oplus \ldots \oplus (0\ldots01_{2i-1}0\ldots01-\frac{3}{2})$$
For $l$ even, $$\bigwedge^{l}\mathbb{\mathbb{V}}^*\otimes
\mathbb{S} \simeq (0\ldots 0-\frac{1}{2}) \oplus (010 \ldots
0-\frac{1}{2})\oplus\ldots \oplus (0\ldots 010-\frac{1}{2}) \oplus(0
\ldots   \ldots 0 \frac{1}{2})$$
$$\oplus (10\ldots  01-\frac{3}{2})\oplus(0010\ldots 01-\frac{3}{2})\oplus \ldots \oplus (0\ldots0101-\frac{3}{2})\oplus(0\ldots 02-\frac{3}{2})$$
For $l$ odd, $$\bigwedge^{l}\mathbb{\mathbb{V}^*}\otimes
\mathbb{S} \simeq (10 \ldots 0-\frac{1}{2}) \oplus (0010 \ldots
0-\frac{1}{2})\oplus\ldots \oplus (0 \ldots  0 1
0-\frac{1}{2})\oplus (0\ldots 0\frac{1}{2})$$
$$\oplus (0\ldots 0 01-\frac{3}{2})\oplus(010\ldots 01-\frac{3}{2})\oplus \ldots \oplus (0\ldots0101-\frac{3}{2})\oplus \ldots (0\ldots 02-\frac{3}{2})$$}

{\it Proof.} Since $\omega: \mathbb{V}\times \mathbb{V} \to
\mathbb{C}$ is a non degenerate $\mathfrak{g}$-invariant bilinear
form, it gives a $\mathfrak{g}$-module isomorphism $\mathbb{V}\simeq
\mathbb{V}^*.$ Thus the decomposition of the product
$\bigwedge^i\mathbb{V}^*\otimes \mathbb{S}$ is equivalent to the
decomposition of $\bigwedge^i\mathbb{V}\otimes \mathbb{S}.$ For
obtaining a further isomorphism, choose a symplectic basis
$\{e_j\}_{j=1}^{2l}$ of $(\mathbb{V},\omega)$ and  define a mapping
$$\phi: \bigwedge^i\mathbb{V}\times \bigwedge^{2l-i}\mathbb{V} \to \mathbb{C}$$ 
 for $i=0,\ldots, 2l$ on homogeneous elements by a formula
$\phi(v_1\wedge \ldots \wedge v_i, w_1\wedge \ldots \wedge
w_{2l-i})=: q \in \mathbb{C},$ if and only if $q e_1\wedge \ldots
\wedge e_{2l}=v_1\wedge \ldots \wedge v_i \wedge w_1 \wedge \ldots
\wedge w_{2l-i},$ where $v_1,\ldots, v_i, w_1,\ldots, w_{2l-i} \in
\mathbb{V}.$ Obviously, one extends the definition by linearity.
 Since the symplectic group $G=Sp(\mathbb{V},\omega)$ is a subgroup
of the special linear group $SL(\mathbb{V})$, we have
$\phi(g v,g w)=g v\wedge g w=\mbox{det}(g)v\wedge w = v \wedge w =\phi(v,w)$ for each $v,w
\in \mathbb{V}$ and $g\in Sp(\mathbb{V},\omega),$ i.e., the mapping
$\phi$ is $Sp(\mathbb{V},\omega)$- and also
$\mathfrak{sp}(\mathbb{V},\omega)$-invariant in the appropriate
manners. Thus $\bigwedge^i \mathbb{V} \simeq
(\bigwedge^{2l-i}\mathbb{V})^*, $ which is naturally isomorphic to
$\bigwedge^{2l-i}\mathbb{V}^*,$ which is in turn isomorphic to
$\bigwedge^{2l-i}\mathbb{V}.$ Thus we need to decompose the spaces
$\bigwedge^i\mathbb{V}\otimes \mathbb{S}$ for $i=0,\ldots, l$ only. After
a straightforward but tedious application of lemma 1 and theorem 2
we would get the decompositions written in the statement of this
lemma. $\Box$

For sake of brevity, let us introduce the following notation. First, let us
define a finite subset $\Xi$ of pairs of non-negative integers.
$$\Xi:=\{(i,j)| i=0,\ldots,l; j=0,\ldots, i\}\cup \{(i,j)|
i=l+1,\ldots,2l, j=0,\ldots, 2l-i\}.$$ Further, let us define
\begin{itemize}
\item[]$\mathbb{E}_{0,2j}:=(0\ldots 0 -\frac{1}{2}),$ $(0,2j)\in \Xi-\{(l,l),(l,l-1)\}$
\item[]$\mathbb{E}_{0,2j+1} := (0\ldots 0 1 -\frac{3}{2}),$ $(0,2j+1)\in \Xi-\{(l,l),(l,l-1)\}$
\item[]$\mathbb{E}_{2i,2j}:=(0 \ldots 0 1_{2j} 0 \ldots 0-\frac{1}{2}),$ $(2i,2j)\in \Xi-\{(l,l),(l,l-1)\}$
\item[]$\mathbb{E}_{2i+1,2j}:=(0 \ldots 0 1_{2j}0\ldots  0 1 -\frac{3}{2}),$ $(2i+1,2j)\in \Xi-\{(l,l),(l,l-1)\}$
\item[]$\mathbb{E}_{2i,2j+1}:=(0 \ldots 0 1_{2j+1} 0 \ldots 0 1-\frac{3}{2}),$ $(2i,2j+1)\in \Xi-\{(l,l),(l,l-1)\}$
\item[]$\mathbb{E}_{2i+1,2j+1}:=(0\ldots 0 1_{2j+1}0\ldots 0-\frac{1}{2} ),$ $(2i+1,2j+1)\in \Xi-\{(l,l),(l,l-1)\}$
\item[]$\mathbb{E}_{l,l-1}:=(0\ldots 0 2 -\frac{3}{2}),$ \, $\mathbb{E}_{l,l}=(0\ldots 0 \frac{1}{2}).$
\end{itemize}

Let us remark, that a little bit more systematic way of defining the
modules $\mathbb{E}_{i,j}$ for $(i,j)\in \Xi$ would be that one, in which the
basis  $\{\epsilon_i\}_{i=1}^l$ is used.

Using this notation, we can reformulate the lemma 2 in the following
way
 $$\bigwedge^{2l-i}\mathbb{V}^*\otimes \mathbb{S}\simeq
 \bigwedge^i\mathbb{V}^*\otimes \mathbb{S} \simeq \bigoplus_{j=0}^i
\mathbb{E}_{i,j}$$ for $i=0,\ldots, l$ or
$$\bigwedge^{i}\mathbb{V}^*\otimes \mathbb{S} \simeq \bigoplus_{(i,j)\in
\Xi}\mathbb{E}_{i,j}$$ for $i=0,\ldots, 2l.$

To visualize the system described by lemma 2, we display a picture for rank $l=3.$
The $i^{th}$  column  corresponds to the space $\bigwedge^{i}\mathbb{V}^*\otimes \mathbb{S}$
and each member of a column corresponds to an irreducible representation in $\bigwedge^{i}\mathbb{V}^*\otimes \mathbb{S}$
with a displayed highest weight.

$$(00-\frac{1}{2})\quad (01-\frac{3}{2})\quad(00-\frac{1}{2})\quad(01-\frac{3}{2})\quad(00-\frac{1}{2})\quad (01-\frac{3}{2})\quad(00-\frac{1}{2})$$
$$(10-\frac{1}{2})\quad(11-\frac{3}{2})\quad(10-\frac{1}{2})\quad(11-\frac{3}{2})\quad(10-\frac{1}{2})$$
$$(01-\frac{1}{2})\quad(02-\frac{3}{2})\quad(01-\frac{1}{2})$$
$$(00\frac{1}{2})$$

In the next theorem,
a decomposition of a tensor product of a higher symplectic spinor module and the 
defining representation $\mathbb{V} \simeq F(\varpi_1)$ over $\mathfrak{sp}(\mathbb{V},\omega) \simeq \mathfrak{sp}(2l,\mathbb{C})$ is described.

{\bf Theorem 3:} Let $\mathfrak{g} \simeq \mathfrak{sp}(2l,\mathbb{C})$ and $\lambda \in A.$ Then
$$L(\lambda) \otimes F(\varpi_1)=\bigoplus_{\mu \in A_{\lambda}}L(\mu),$$
where $A_{\lambda}:=A \cap \{\lambda + \nu| \nu \in \Pi(\varpi_1)
\}$ and $\Pi(\varpi_1)$ is the saturated set of weights of the
defining representation.\footnote{One can easily compute,that
$\Pi(\varpi_1)=\{\pm\epsilon_i| i=1,\ldots,l\}.$}

{\it Proof.} See Kr\'{y}sl, \cite{Krysl} or \cite{KryslLie}. $\Box$

Let us remark, that the proof of this theorem is based on the so called Kac-Wakimoto formal character formula of Kac and Wakimoto published in \cite{KW} and on some results of Humphreys, see \cite{Hum2}, who specified results of Kostant from \cite{K} on tensor products of finite and infinite dimensional modules admitting a central character.

In the next lemma, a property is formulated, which is valid for an arbitrary simple Lie algebra $\mathfrak{g}.$

{\bf Lemma 3:} Let $X$ be a $\mathfrak{g}$-module and $\mathbb{V},
\mathbb{W} \subseteq X$ its two $\mathfrak{g}$-submodules (of finite
or infinite dimension). Suppose $\mathbb{V}=\mathbb{V}_{1}\oplus
\ldots \oplus \mathbb{V}_{a}$ and
$\mathbb{W}=\mathbb{W}_1\oplus\ldots\oplus \mathbb{W}_{b}$ are
decompositions into irreducible $\mathfrak{g}$-submodules. Define a
subset $I$ of the set $\{1,\ldots, a\}$ by the prescription
$I:=\{i\in \{1,\ldots,a\}| \exists j \in \{1,\ldots,
b\}:\mathbb{V}_i \simeq \mathbb{W}_j\}.$ Then

$$\mathbb{V} \cap \mathbb{W}\subseteq \bigoplus_{i\in I}\mathbb{V}_i.$$

{\it Proof.} Let us consider the projections $p_i: \mathbb{V} \to
\mathbb{V}_i,$ $i=1,\ldots, a$ and $q_j:\mathbb{W}\to \mathbb{W}_j,$
$j=1,\ldots, b.$ Suppose, that we have defined the projections $p_i$
also on the space $\mathbb{W}$  in the following way. One can easily show, that a finite
direct sum of $\mathfrak{g}$-modules is actually completely
reducible (see Kr\'{y}sl \cite{Krysl}). Thus there is a (generally
non-unique) $\mathfrak{g}$-submodule $\mathbb{U},$ such that
$(\mathbb{V}\cap \mathbb{W})\oplus\mathbb{U} = \mathbb{W}$ is a direct
sum of $\mathfrak{g}$-modules. Therefore given any $x\in
 \mathbb{W},$ we can write it as a sum $x=v+u,$ where
$v\in \mathbb{V}\cap \mathbb{W}$ and $u \in \mathbb{U},$ in a unique way and define
$p_i(x):=p_i(v).$ Now, take an element $x\in \mathbb{V}\cap
\mathbb{W}.$ We have $x=\sum_{i=1}^a\sum_{j=1}^b p_i(q_{j
|\mathbb{W}_j}(x)).$
To get a contradiction, suppose there are  elements $i \notin I$ and $j\in \{1,\ldots,b\}$ such
that $p_i(q_j(x))\neq 0.$   Thus we have a $\mathfrak{g}$-module
homomorphism $R:=p_i\circ q_{j}|_{\mathbb{W}_j}: \mathbb{W}_j \to \mathbb{V}_i,$ which is nonzero. The
classical Schur lemma type argument shows that $R$ is an isomorphism
of $\mathbb{W}_j$ and $\mathbb{V}_i,$ which contradicts the condition $i \notin I.$ $\Box$

{\bf Remark:} The proof of the above written lemma does not use any
information about the Lie algebra over which we took the module $X,$
thus it could be generalized to each module over a general algebraic
structure (group, commutative, associative,   super-Lie algebra
e.t.c.), which admits modules over itself.
   Let us note, that the
statement of the lemma can be improved in an easy way. To see the
weakness of the lemma, consider two nonzero equivalent
representations $\mathbb{V}_1$ and $\mathbb{V}_2,$ form their direct sum $\mathbb{V}=\mathbb{V}_1\oplus
\mathbb{V}_2$ and suppose a submodule $\mathbb{W}\simeq \mathbb{V}_1\simeq \mathbb{V}_2,$ for which $\mathbb{V}_1 \neq \mathbb{W} \neq \mathbb{V}_2,$
 is given.
  Then clearly $\mathbb{V}\cap \mathbb{W} = \mathbb{W} \subsetneq \mathbb{V}_1\oplus \mathbb{V}_2,$ but the lemma
gives only $\mathbb{V}\cap \mathbb{W} \subseteq \mathbb{V}_1\oplus \mathbb{V}_2.$ We will not try to improve
this lemma, because first we shall need it only in the above written
form and second the reformulation would be a bit inefficient because
of its increased length.

\section{Spin symplectic geometry}

We shall begin with a short observation about covariant derivatives on vector bundle valued forms and then we are going to 
consider basic aspects of  metaplectic structures.

{\bf Lemma 4:} Let $p:F \to M$  be a smooth vector bundle equipped
by a vector bundle connection $\nabla^F$.  Consider a subbundle $E
\subseteq \bigwedge^i T^*M \otimes F \to M$ for some $i=0,\ldots, \mbox{dim}M$ and a section $s \in
\Gamma(M,E).$ Since $s$ could be viewed as an exterior differential
form with values in $F,$ the covariant derivative $d^{\nabla^F}$
could be applied. Then $$d^{\nabla^F}s \in \Gamma(M,(T^*M\otimes  E)
\cap (\bigwedge^{i+1}T^*M\otimes F)).$$

 {\it Proof.}   This is an easy observation. We know that $s \in \Gamma(M,   \bigwedge^iT^*M \otimes F)$ and therefore
$d^{\nabla^F}s \in \Gamma(M,\bigwedge^{i+1}T^*M \otimes F).$ The assumption $s \in \Gamma(M, E)$ implies $d^{\nabla^F}s \in \Gamma(M,T^*M\otimes E).$
Summing up, we obtain $d^{\nabla^F}s \in \Gamma(M,(T^*M\otimes  E) \cap (\bigwedge^{i+1}T^*M\otimes F)).$ $\Box$

To define a metaplectic structure, we will use a definition of
Katharina Habermann from \cite{H}, which is quite analogous to the
Riemannian case. Now, let $(\mathbb{V}_0,\omega)$ be a real
symplectic space of dimension $2l.$ Let $\tilde{G}_0$ be a nontrivial 2-fold
covering of the group $G_0=Sp(\mathbb{\mathbb{V}}_0,\omega) \simeq Sp(2l,\mathbb{R}),$   thus
$\tilde{G}_0\simeq Mp(2l,\mathbb{R})$ (the metaplectic group) and
fix a 2:1 covering $\lambda: \tilde{G}_0 \to G_0.$

{\bf  Definition 2:} Let $(M,\omega)$ be a   symplectic manifold of
 dimension $2l.$ Let $(p_1: \mathcal{P} \to M, Sp(2l,\mathbb{R}))$
be a principal fiber bundle of symplectic repers (in $TM$) and
$(p_2:\mathcal{Q} \to M, Mp(2l,\mathbb{R}))$ be a principal fiber
bundle with a structure group $Mp(2l,\mathbb{R}).$ We call a
surjective bundle homomorphism $\Lambda: \mathcal{Q} \to
\mathcal{P}$ (over the identity on $M$) {\it metaplectic structure}, if
the following diagram commutes.

$$\begin{xy}\xymatrix{
Mp(2l,\mathbb{R}) \times \mathcal{Q} \ar[dd]^{\lambda\times \Lambda} \ar[r]&   \mathcal{Q} \ar[dd]^{\Lambda} \ar[dr]^{p_2} &\\
                                                            & &M\\
Sp(2l,\mathbb{R}) \times \mathcal{P} \ar[r]   & \mathcal{P} \ar[ur]_{p_1}
}\end{xy}$$

Let $(M, \omega)$ be a symplectic manifold of dimension $2l.$ It is
well known that there is no unique {\it symplectic connection}.
Symplectic connection is a torsion-free  connection $\nabla$ on the
tangent bundle which preserves the symplectic structure $\omega,$
i.e., $$\omega(\nabla_XY,Z)+\omega(Y,\nabla_{X}Z)=X\omega(Y,Z)$$ for each $X,
Y, Z \in \Gamma(M,TM),$ see, e.g., Habermann \cite{H}.  Nevertheless,
 we may take any symplectic connection $\nabla$ and associate to it a
principal bundle connection $Z^{\mathcal{P}}:T\mathcal{P} \to
\mathfrak{g}_0.$ For this connection, there is a lifted connection
$Z^{\mathcal{Q}}: T\mathcal{Q}\to \tilde{\mathfrak{g}_0}\simeq
\mathfrak{g}_0$ on the metaplectic structure. To this connection,
$Z^{\mathcal{Q}}$ we can associate a linear connection $\nabla^S$ on
the spinor bundle $p:S\to M.$ By the spinor bundle $p:S\to M$ we
mean the associated vector bundle
$S:=\mathcal{Q}\times_{\tilde{G}_0}\mathbb{S}$ to the principle
$\tilde{G}_0$-bundle via the basic symplectic spinor representation
$\mathbb{S}.$ In this way we may construct the covariant
derivative $d^{\nabla^S}.$ (For the correctness of this definition,
see Kashiwara, Vergne \cite{KV}, where certain globalization of the basic spinor modules to the group $\tilde{G}_0$ is described.)

For the sake of brevity, let us introduce the notation for
the following associated vector bundles
$E_{i,j}:=\mathcal{Q}\times_{\tilde{G}_0}\mathbb{E}_{i,j}$ for all
$(i,j)\in \Xi.$ For technical reasons, we define
$E_{i,j}:=0$ if $(i,j) \notin \Xi.$

{\bf Theorem 4:} Let $\nabla$ be a torsion-free symplectic connection on the tangent bundle of the (symplectic) base  manifold M of a metaplectic structure. Let us denote the induced covariant derivative for the associated symplectic spinor bundle $p: S \to M$ by $d^{\nabla^{S}}.$
 Then with the notation introduced above,
$$d^{\nabla^S}:\Gamma(M,E_{i,j})\to \Gamma(M,E_{i+1,j-1}\oplus E_{i+1,j}\oplus E_{i+1,j+1})$$
 for all $(i,j) \in \Xi.$

{\it Proof.} We must prove that the image is of the form described
by the statement of the theorem.  Lemma 4 (for $F:=S, E:=E_{i,j} \subseteq
\bigwedge^{i}T^*M \otimes S$ and a connection $\nabla^S$ induced by
a symplectic connection $\nabla$ on $TM$) implies, that
$d^{\nabla^S}s \in \Gamma(M,(E_{i,j}\otimes T^*M) \cap
\bigwedge^{i+1}T^*M\otimes S)$ for a section $s \in
\Gamma(M,E_{i,j}).$ The vector bundle $T^*M\otimes E_{i,j} \cap
\bigwedge^{i+1}T^*M\otimes S$ is isomorphic to the associated vector
bundle $\mathcal{Q}\times_{\tilde{G}}(\mathbb{\mathbb{V}}^*\otimes
\mathbb{E}_{i,j}\cap \bigwedge^{i+1}\mathbb{V}^*\otimes
\mathbb{S}).$ Our strategy is to use lemma 3 to compute the intersection. Therefore we need
to find the decomposition of the modules $\mathbb{\mathbb{V}}^*\otimes
\mathbb{E}_{i,j}$ and $\bigwedge^{i+1}\mathbb{V}^*\otimes
\mathbb{S}$ into irreducible summands. The latter was done in lemma
2. We shall use  theorem 3 to decompose $\mathbb{V}^*\otimes \mathbb{E}_{i,j}
\simeq \mathbb{V}\otimes \mathbb{E}_{i,j}.$ There are in
principle 6 forms of $\mathbb{E}_{i,j}:$ $L(0\ldots 0
-\frac{1}{2}),$ $L(0\ldots 01-\frac{3}{2}),$ $L(0 \ldots 0 1_k 0
\ldots 0 -\frac{1}{2}),$ $L(0\ldots 0 1_k 0 \ldots 0
1-\frac{3}{2}),$ $L(0\ldots 0 2-\frac{3}{2})$ and $L(0 \ldots 0
\frac{1}{2}).$ Because we will be more careful and will distinguish between
odd and even subscripts $i,j$ for the space $\mathbb{E}_{i,j}$ at the beginning of our analysis,
 we will be investigating eleven cases actually.

\begin{itemize}

\item[1.)] For $\mathbb{E}_{2i,0}=(0\ldots 0  -\frac{1}{2})$
we obtain that $\mathbb{E}_{2i,0}\otimes \mathbb{V} = (10\ldots 0-\frac{1}{2}) \oplus (0\ldots 0 1 -\frac{3}{2})$

 \begin{itemize}
 \item For $i=0,\ldots, l-1,$ we obtain that $((10\ldots 0 - \frac{1}{2})\oplus (0\ldots 0 1 -\frac{3}{2}))\cap \bigwedge^{2i+1}\mathbb{V}\otimes           \mathbb{S}  \subseteq  (10\ldots 0 - \frac{1}{2})\oplus (0\ldots 0 1 -\frac{3}{2})=\mathbb{E}_{2i+1,-1}\oplus \mathbb{E}_{2i+1,0}\oplus  \mathbb{E}_{2i+1,1}$  (the first summand is zero by definition).
 \item For $i=l,$ we obtain that the intersection is zero, because $\bigwedge^{2i+1}\mathbb{V}\otimes \mathbb{S}$ is zero.
 But in this case the vector space  $\mathbb{E}_{2l,-1}\oplus \mathbb{E}_{2l+1,0}\oplus \mathbb{E}_{2l+1,1}$ is zero, too (by definition).
 \end{itemize}

In each cases, we have obtain that $\mathbb{E}_{2i,0}\otimes \mathbb{V}^* \subseteq \mathbb{E}_{2i,-1}
\oplus\mathbb{E}_{2i,0}\oplus\mathbb{E}_{2i,1}$  according to the statement.

\item[2.)] For $\mathbb{E}_{2i+1,0}=(0\ldots 0 1 -\frac{3}{2})$ we obtain
that $\mathbb{E}_{2i+1,0}\otimes \mathbb{V} = (10\ldots 0 1 -\frac{3}{2})\oplus (0\ldots 0-\frac{1}{2})\oplus (0\ldots 010-\frac{3}{2}).$

\begin{itemize}
\item For $i=0,\ldots,l-2$ we get $((10\ldots 0 1 -\frac{3}{2})\oplus (0\ldots 0-\frac{1}{2})\oplus (0\ldots 010-\frac{3}{2}))\cap \bigwedge^{2i+1}\mathbb{V}\otimes \mathbb{S}=(0\ldots 0-\frac{3}{2})\oplus (10\ldots 0 1 -\frac{3}{2})=\mathbb{E}_{2i+2,-1}\oplus \mathbb{E}_{2i+2,0} \oplus \mathbb{E}_{2i+2,1},$ because the first summand is zero (by definition).
\item For $i=l-1$ we obtain $((10\ldots 0 1 -\frac{3}{2})\oplus (0\ldots 0-\frac{1}{2})\oplus (0\ldots 010-\frac{3}{2}))\cap \bigwedge^{2l}\mathbb{V}\otimes \mathbb{S}=(0\ldots 0 -\frac{1}{2})=\mathbb{E}_{2l,-1}\oplus \mathbb{E}_{2l,0}\oplus \mathbb{E}_{2l,1},$ because
the first and last summands are zero (by definition).
\end{itemize}

\item[3.)] For  $\mathbb{E}_{2i,2j}=(0\ldots 0 1_{2j}
 0 \ldots 0 -\frac{1}{2})$ ($l-2\geq 2j>0$ is to be supposed, $j=0$ has been already handled)
 we obtain  $\mathbb{E}_{2i,2j}\otimes \mathbb{V}=
 (0\ldots 0 1_{2j-1} 0 \ldots 0-\frac{1}{2})\oplus (0\ldots 0 1_{2j+1}0\ldots 0 -\frac{1}{2})\oplus
(0\ldots 01_{2j}0\ldots 1-\frac{3}{2})\oplus (1 0\ldots 0 1_{2j}0\ldots 0 -\frac{1}{2}).$

\begin{itemize}
\item For $2i+2j<2l$, the intersection $((0\ldots 0 1_{2j-1} 0 \ldots 0-\frac{1}{2})
\oplus$ \\ $(0\ldots 0 1_{2j+1}0\ldots 0 -\frac{1}{2})\oplus
(0\ldots 01_{2j}0\ldots 1-\frac{3}{2})\oplus (1 0\ldots 0
1_{2j}0\ldots 0 -\frac{1}{2}))$ $\cap
\bigwedge^{2i+1}\mathbb{V}\otimes \mathbb{S}=
 (0\ldots 0 1_{2j-1} 0 \ldots 0-\frac{1}{2})\oplus (0\ldots 0 1_{2j+1}0\ldots 0 -\frac{1}{2})\oplus
(0\ldots 01_{2j}0\ldots 1-\frac{3}{2}) = \mathbb{E}_{2i+1,2j-1}\oplus \mathbb{E}_{2i+1,2j}\oplus \mathbb{E}_{2i+1,2j+1}.$
\item  For $2i+2l=2l,$ we get that the intersection equals  $(0\ldots 0 1_{2j-1} 0 \ldots 0-\frac{1}{2}).$
\end{itemize}

In each case, we have obtained that the intersection is in $\mathbb{E}_{2i+1,2j-1}\oplus E_{2i+1,2j}\oplus \mathbb{E}_{2i+1,2j+1}.$
\end{itemize}

The remaining cases will be handled not so carefully. We will write only the result of the appropriate decomposition and not the result of the intersection. In each cases, one can compute the intersection like in the previous ones and check that the condition for the intersection to be a vector subspace of $\mathbb{E}_{i+1,j-1}\oplus\mathbb{E}_{i+1,j}\oplus \mathbb{E}_{i+1,j+1}$ for appropriate $i,j$ is fulfilled. We will also not distinguish between the parities of $i,j$

\begin{itemize}
\item[4.)] For $\mathbb{E}=(0 \ldots 0 1_k 0\ldots 0 1 -\frac{3}{2})$ ($k>1$), we get $\mathbb{E}\otimes \mathbb{V}=
(0\ldots 0 1_{k-1} 0 \ldots 0 1-\frac{3}{2}) \oplus (0\ldots 0 1_{k+1} 0 \ldots 0 1 -\frac{3}{2}) \oplus (1 0 \ldots 0 1_k 0 \ldots 0 1-\frac{3}{2}) \oplus (0\ldots 0 1_k 0 \ldots 0 -\frac{1}{2})$

\item[5.)] For $\mathbb{E}=(0 \ldots 0 1_k 0\ldots 0 -\frac{1}{2})$ ($k>1$), we get $\mathbb{E}\otimes \mathbb{V}=
(0\ldots 0 1_{k-1} 0 \ldots 0 -\frac{1}{2}) \oplus (0\ldots 0 1_{k+1} 0 \ldots 0  -\frac{1}{2}) \oplus (1 0 \ldots 0 1_k 0 \ldots 0 -\frac{1}{2}) \oplus (0\ldots 0 1_k 0 \ldots 0 1-\frac{3}{2})$

\item[6.)] For $\mathbb{E}=(10\ldots 0 -\frac{1}{2})$ we get $\mathbb{E}\otimes \mathbb{V}=(0\ldots 0-\frac{1}{2})\oplus (010\ldots 0-\frac{1}{2})\oplus
(1\ldots 0 1-\frac{3}{2}).$

\item[7.)] For $\mathbb{E}= (10\ldots 0 1-\frac{3}{2})$ we get $\mathbb{E}\otimes \mathbb{V}=(0\ldots 01-\frac{3}{2})\oplus (010\ldots 01-\frac{3}{2})\oplus (1\ldots 0 -\frac{1}{2}).$

\item[8.)] For $\mathbb{E}=(0\ldots 0 1 -\frac{1}{2})$ we get $\mathbb{E}\otimes \mathbb{V}=(1 0 \ldots 0 1 -\frac{1}{2}) \oplus (0 \ldots 0 2 -\frac{3}{2}) \oplus (0 \ldots 0 \frac{1}{2}).$

\item[9.)]For $\mathbb{E}=(0\ldots 0 1 1 -\frac{3}{2})$ we get
$\mathbb{E}\otimes \mathbb{V}=(1 0 \ldots 0 1 1 -\frac{3}{2}) \oplus (0 \ldots 0 2 -\frac{3}{2}) \oplus (0 \ldots 0 1 0 -\frac{1}{2}).$

\item[10.)] For $\mathbb{E}=(0\ldots 0 2 -\frac{3}{2})$ we get $\mathbb{E}\otimes \mathbb{V}=
(10\ldots 0 2-\frac{3}{2}) \oplus (0\ldots 0 1 1 -\frac{3}{2})\oplus (0 \ldots 0 3 -\frac{5}{2})\oplus (0\ldots 0 1 -\frac{1}{2}).$

\item[11.)] For $\mathbb{E}=(0\ldots 0 \frac{1}{2})$ we get $\mathbb{E}\otimes \mathbb{V} = (10\ldots 0 \frac{1}{2})\oplus (0 \ldots 0  -\frac{1}{2}).$

\end{itemize}

The irreducible representation of $\mathfrak{g}$ could be define
also for the split real form
$\mathfrak{g}_0=\mathfrak{sp}(2l,\mathbb{R}).$ The irreducibility
does not change and also the decompositions remain the same. Let us
remark, that according to the result of Kashiwara, Vergne \cite{KV},
there are some $L^2$-globalizations of our representation to the
metaplectic group $\tilde{G}_0$ and that there are also the canonically defined ones. The decompositions do not change if
we turn our attention to the $(\mathfrak{g}_0,\tilde{K})$-structure,\footnote{$\tilde{K}$ 
is the maximal compact subgroup of
$\tilde{G_0}\simeq Mp(2l,\mathbb{R}),$ i.e., the $\lambda$-preimage of the maximal
compact subgroup $K$ of $Sp(2l,\mathbb{R}),$ which is isomorphic to the unitary group, $K \simeq U(2l).$}
because $\tilde{K}$ is connected (see Baldoni \cite{Baldoni}), and
do not change even if we take the appropriate globalization (e.g., the Casselman-Wallach, i.e., the minimal one). $\Box$

In  the next picture, the system described by  theorem 4 is displayed for $l=3,$ i.e., for metaplectic structures over a six dimensional symplectic
manifold $(M^6,\omega).$

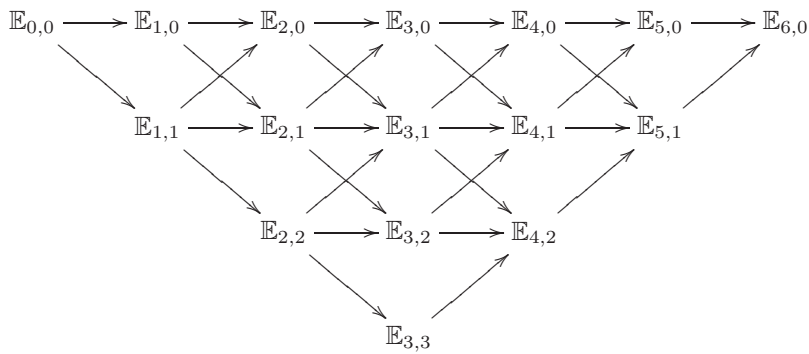
$$\xymatrix{\mathbb{E}_{0,0}\ar[r]\ar[dr] &\mathbb{E}_{1,0}\ar[r]\ar[dr] &\mathbb{E}_{2,0}\ar[r]\ar[dr] &\mathbb{E}_{3,0} \ar[r]\ar[dr] &\mathbb{E}_{4,0}\ar[r]\ar[dr] &\mathbb{E}_{5,0}\ar[r] &\mathbb{E}_{6,0}\\
&\mathbb{E}_{1,1} \ar[ur]\ar[r]\ar[dr] &\mathbb{E}_{2,1}\ar[r]\ar[dr] \ar[ur]            &  \mathbb{E}_{3,1}\ar[r]\ar[dr]\ar[ur]  & \mathbb{E}_{4,1} \ar[r]\ar[ur]   &\mathbb{E}_{5,1}\ar[ur] & \\
 &               & \mathbb{E}_{2,2} \ar[r] \ar[ur] \ar[dr] & \mathbb{E}_{3,2} \ar[ur] \ar[r] & \mathbb{E}_{4,2}\ar[ur]&&\\
&&& \mathbb{E}_{3,3}\ar[ur]&&&}
$$

Having the analogous result in the Riemannian case in mind, we are entitled to
call the horizontally going operators in the first arrow {\it symplectic Dirac
operators}, that ones going from the first arrow down-right {\it symplectic
twistor} and  in the second arrow the horizontally going ones {\it symplectic
Rarita-Schwinger operators}. The horizontally going operators on the remaining arrows could be eventually
called {\it symplectic generalized Rarita-Schwinger operators.}

The further research could be devoted to other  real symplectic
groups, various types of globalizations  of the modules in question,
to a coordinate-way description of the operators, we have obtained,
and to their analytic properties.

 \thanks{I am very grateful to Vladim\'ir Sou\v{c}ek for
 motivations  coming  from Riemannian geometry
 and orthogonal spin structures.
 The author is also grateful to the Grant Agency of Czech Republic for the support from the grant
for young researchers GA\v{C}R 201/06/P223.}

\end{document}